\documentclass[a4paper, 11pt]{amsart}

\usepackage{dsfont}
\usepackage{amsmath, amsthm, amssymb}
\usepackage{fullpage}
\usepackage{xypic}
\usepackage{graphicx}
\usepackage{color}
\usepackage{mathtools}
\usepackage[usenames,dvipsnames]{xcolor}
\usepackage[latin1]{inputenc}
\usepackage{indentfirst}
\usepackage{epstopdf}
\usepackage{subfigure}
\usepackage{emptypage}
\usepackage{hyperref}
\hypersetup{
    colorlinks=true, %set true if you want colored links
    linktoc=all,     %set to all if you want both sections and subsections linked
    linkcolor=blue,  %choose some color if you want links to stand out
    citecolor=JungleGreen,
}

\newtheorem{theorem}{Theorem}%[section]

\newtheorem{theorem*}{Theorem}
\newtheorem{question*}[theorem*]{Question}
\newtheorem{conjecture*}[theorem*]{Conjecture}
\newtheorem{corollary*}[theorem*]{Corollary}
\newtheorem{theorem*e}{Teorema}
\newtheorem{question*e}[theorem*e]{Pregunta}
\newtheorem{conjecture*e}[theorem*e]{Conjetura}
\newtheorem{corollary*e}[theorem*e]{Corolario}

\renewcommand{\det}{\mathrm{det}}
\newcommand{\id}{\mathrm{id}}
\newcommand{\dist}{\mathrm{dist}}

\begin{document}

\title{Parrondo's paradox for homeomorphisms}
%\title{Parrondo's paradox in planar dynamics}
\author{A. Gasull}
\address{
Departament de Matem\`atiques, Universitat Aut\`onoma de Barcelona and Centre de Recerca Matem\`atica \\
Campus de Bellaterra
08193 Bellaterra, Barcelona, Spain.
}
\email{gasull@mat.uab.cat}

\author{L. Hern\'andez-Corbato}
\address{Departamento de \'Algebra, Geometr\'{\i}a y Topolog\'{\i}a\\ Universidad Complutense de Madrid and Instituto de Ciencias Matematicas CSIC--UAM--UCM--UC3M \\ Madrid \\ Spain.}

\email{luishcorbato@mat.ucm.es}

\author{F. R. Ruiz del Portal}
\address{Departamento de \'Algebra, Geometr\'{\i}a y Topolog\'{\i}a\\ Universidad Complutense de Madrid \\ 28040 Madrid \\ Spain.}
\email{rrportal@ucm.es}

\begin{abstract}
We construct two planar homeomorphisms $f$ and $g$ for which the
origin is a globally asymptotically stable fixed point whereas for $f \circ g$
and $g \circ f$ the origin is a global repeller.
Furthermore, the origin remains a global repeller for
the iterated function system generated by $f$ and $g$ where each of
the maps appears with a certain probability. This planar
construction is also extended to any dimension greater
than 2 and proves for first time the appearance of the
Parrondo's dynamical paradox in odd dimensions.
\end{abstract}

\maketitle



\section{Introduction and main results}

The  \emph{Parrondo's paradox} is a well-known paradox in game
theory, that in a few words affirms that {\it a combination of
losing strategies can become a winning strategy,} see~\cite{HA,Par}.
In the dynamical context,  when we study the stability of fixed
points, the role of being a winning or a losing strategy can be
replaced by being an attracting or repelling fixed point.
A word of caution, throughout this note we use the term attracting (or attractor) and repelling (or repeller) as a synonym of asymptotically stable for a map and its inverse, respectively.
Hence,
for a fixed class of maps $\mathcal{C},$ from $\mathbb{R}^k$ into
itself, we will say that a pair of maps $f,g\in\mathcal{C}$
exhibit a \emph{dynamical Parrondo's paradox}
if they have a common fixed point at which the maps are locally invertible and the fixed point is locally asymptotically stable for $f$ and $g$ but it is a repeller for the
composite maps $g\circ f$ and $f\circ g.$ Notice that, 
$g\circ f$ and $f\circ g$ are conjugate near the fixed
point because, locally,  $f\circ g= g^{-1}\circ g\circ f \circ g.$

As shown in \cite{CGM18}, the dynamical Parrondo's paradox can arise 
%In \cite{CGM18} it is proved that  a dynamical Parrondo's paradox happens
when $k$  is even and $\mathcal C$ is the class of
polynomial maps. On the contrary, it is also proved in \cite{CGM18} that the paradox does not appear when $k=1$ and $\mathcal C$ is the class of
analytic maps. Notice also that the paradox is also impossible
for any $k$ when $\mathcal C$ is  the class of maps for which the
common fixed point is hyperbolic. Indeed,
 given two $k\times k$ matrices with all their
    eigenvalues with modulus smaller than 1, $A$ and $B$, it holds
that $|\det(A)|<1,|\det(B)|<1  $ and hence $|\det(A B)|<1$. As a
consequence, when two maps share  a common fixed point ${\bf x}$
which is asymptotically stable for both of them and the maps are of class
$\mathcal{C}^1$ at   ${\bf x},$ then, generically, for $g\circ f$
and $f\circ g$ the fixed point~${\bf x}$ is, in both cases, either
locally asymptotically stable or of saddle type, but it can never be
repeller. Examples of saddle type points for $g\circ f,$ when $f$
and $g$ are linear maps, are given in \cite{BTT} and \cite[p. 8]{J}.

For the sake of completeness, and to compare it with our result, we recall the
example in \cite[Ex. 7]{CGM18} for $k=2,$
    \begin{align*}
    f(x,y)&=\big(-y+2{x}^{2}+6xy,x-3{x}^{2}+2xy+3{y}^{2}\big),\\
    g(x,y)&=\big({x}/{2}-{\sqrt{3}}y/{2}-x\,(x^2+y^2),
    {\sqrt{3}}x/{2}+{y}/{2}-y\,(x^2+y^2)\big).
    \end{align*}
It can be proved that the origin is a locally asymptotic stable  fixed point for
    $f$ and $g$ and  the origin is a repelling fixed    point for  $g\circ f$ by
      computing the so called Birkhoff stability constants for the three maps. Notice that the
      dynamics near the fixed points is of rotation type.
Taking the product of these maps $j$ times with themselves
%Taking $j$ couples of maps of this type
we trivially obtain examples of pairs of maps exhibiting the dynamical Parrondo's paradox for all $k=2j.$

It is worth to mention that in  \cite{CLP} a different type of
dynamical Parrondo's paradox is considered. The authors combine
periodically two 1-dimensional maps $f$ and $g$ to  give rise to
chaos or order.

The main goal of this paper is to give examples  of the dynamical
Parrondo's paradox when $\mathcal C$ is the class  of homeomorphisms
and to fill the the lack of examples in odd dimension. We
prove that for all $k\ge2$ there are pairs of maps that realize the {dynamical Parrondo's paradox. While the approach of \cite{CGM18} is mainly
analytic, our point of view is more qualitative. Moreover, the behaviour of our maps near the fixed point is not of
rotation type and there does not seem to be a clear path to make them smooth or analytic.

 \begin{theorem}\label{th:main1}
    For any $k\ge2$, there are pairs of homeomorphisms from $\mathbb{R}^k$ to itself that exhibit the dynamical Parrondo's
    paradox. 
However, for $k=1$ the paradox never arises.
%    {\color{blue} Moreover, the Parrondo's paradox does not
%    happen for $k=1$ in the class of homeomorphisms.}
 \end{theorem}

It is interesting to remark that Theorem \ref{th:main1}, and in
consequence, the dynamical Parrondo's paradox is also relevant from
the point of view of 2-periodic discrete dynamical systems. In
particular, these systems are  good models for describing the
dynamics of biological systems under periodic fluctuations  due
either to external disturbances or effects of seasonality, see
\cite{ES1,ES2,FS,SR1,SR2} and the references therein.

As a byproduct of our construction of the 2-dimensional example of
dynamical Parrondo's paradox we will also prove that, almost surely,
every orbit of the iterated function system generated by $f$ and $g$
is repelled from the origin, where $f$ and $g$ are
essentially the maps constructed in Theorem \ref{th:main1} and appear with
certain respective probabilities $p$ and $1-p$. The result carries onto higher dimensions as well.  To be more precise,
consider the space $\{0, 1\}^{\mathbb N}$ equipped with the
probability measure $\mu$ defined as the product of the Bernoulli
probability measures, $\mu_B,$ in each factor. Recall that for the
Bernoulli distribution $\operatorname{B}(p),$  $\mu_B(1)=p$ and
$\mu_B(0)=q=1-p,$ for some $p\in[0,1].$ We prove:

 \begin{theorem}\label{th:main2}
    For any $k \ge 2$ and $0<p<1,$ there exist homeomorphisms $f_0$ and $f_1$ from
    $\mathbb{R}^k$ into itself such that:
\begin{itemize}
     \item [(i)]The origin is fixed and globally asymptotically stable for $f_0$ and $f_1$.

     \item [(ii)] For ${\bf 0}\ne {\bf x}\in\mathbb{R}^k$ and  for $\mu$-almost all $(a_n) \in \{0, 1\}^{\mathbb N}$ the
orbit $\{F_{a_n, \cdots, a_1, a_0}({\bf x})\}_{n \ge 1}$  is
repelled from the origin,  where $F_{a_n,\cdots ,a_1 , a_0} =
f_{a_n} \circ \ldots \circ f_{a_1} \circ f_{a_0},$ for $n \ge 1,$
$f_0$ appears with probability $p$ and $f_1$ with probability $1-p.$
\end{itemize}
 \end{theorem}

As we will see in the proof, for any $0<p<1,$ each homomorphism
$f_0$ and $f_1$ will have  a region where the radial component of
the points increases and another one where it decreases.
The largest of these variations corresponds to the increasing region, which is in turn considerably bigger in size than the decreasing region.
%Moreover, between these variations the larger  is the one corresponding to the increasing part. 
Their difference becomes larger and tends to infinity when $p$ approaches 0 or 1.

\section{Definition of $f$ and $g$ and proof of Theorem \ref{th:main1}}

We will split the proof in three cases: $k=1,$ $k=2$ and $k>2.$

\subsection{Proof of Theorem \ref{th:main1} for $k=1$}

Let us proceed by contradiction. Suppose that $f$ and $g$ are homeomorphisms of $\mathbb R$, $0$ is a locally attracting fixed point for $f$ and $g$ and a locally repelling fixed point for $f \circ g$ and $g \circ f$. Assume further, for simplicity, that $f$ and $g$ reverse orientation, the other cases are handled similarly. Then:
\begin{itemize}
\item[$(i)$] $g$ is monotone decreasing, so if $y < x < 0$ then $0 < g(x) < g(y)$.
\item[$(ii)$] Since $0$ is locally attracting for $f$ and $g$, for any $y < 0 < x$ close to $0$ we have that 

\centerline{$y < f \circ f(y) < 0 < g \circ g(x) < x$.}
\item[$(iii)$] Since $0$ is locally repelling for $f \circ g$ and $g \circ f$, for any $x,y > 0$ close to $0$ we have that
\centerline{$x < f \circ g(x)$ \enskip and \enskip $y < g \circ f(y)$.}
\end{itemize}
These properties put together yield a contradiction because for small positive $u > 0$:

\centerline{
$u < f \circ g(u) < g \circ f(f\circ g(u)) = g(f \circ f(g(u)) < g(g(u)) < u$
}
\noindent
where the first two inequalities are consequence of $(iii)$ and the last two inequalities are consequence $(i)$ and $(ii)$.

However, notice that it is possible to construct an example in which the origin is semistable for $f \circ g$ (and also for $g \circ f$) while it is asymptotically stable for $f$ and $g$:
\[
f(x) = \begin{cases} -2x & \text{if }x \le 0, \\ -x/3 &\text{if } x \ge 0, \end{cases}
\enskip \enskip\enskip\enskip
g(x) = \begin{cases} -x/3 & \text{if }x \le 0, \\ -2x &\text{if } x \ge 0. \end{cases}
\enskip \enskip\Rightarrow\enskip\enskip
f \circ g(x) = \begin{cases} x/9 & \text{if }x \le 0, \\ 4x &\text{if } x \ge 0. \end{cases}
\]

\subsection{Proof of Theorem \ref{th:main1} for $k=2$}\label{ss:fg}

In our example the maps $f$ and $g$ are conjugate. We first focus on the definition of $f$, expressed in polar coordinates.
%The map $f$ will be defined using polar coordinates in the plane.
The first elements of a typical orbit under $f$ will drift away from the
origin (the radial coordinate increases) until it reaches a trapping
sector in which the orbit remains forever and is slowly attracted to
the origin (the radial coordinate steadily decreases). The dynamics
of the angular coordinate is independent from the values of the radial
coordinate and forces every orbit to be eventually contained in the
trapping sector. Note that, globally, $f$ looks mostly
expanding because the size of this sector and the speed of
convergence to the origin therein are relatively small.

%Despite the size of this sector and the speed of convergence to the origin therein are
%small, the dynamics in the angular coordinate forces every orbit to be eventually contained in the trapping sector where it slowly converges to the origin.

Let us write out the details. Identify $\mathbb R^2 \setminus \{\textbf{0}\}$
with  the cylinder $\mathbb R \times \mathbb R /\mathbb Z$ so as to
use polar coordinates $(r, \theta)$ where $r \in \mathbb R$. The
origin corresponds to the lower end ($r = -\infty$) of the cylinder.
Notice that this is not the typical range for the radial coordinate
but it will later ease our computations. Let $I \subset \mathbb R /
\mathbb Z$ be an interval centered at $\bar{0}$ (here $\bar{0}$ is
  used to denote the neutral element in $\mathbb R / \mathbb Z$) and such that $I
\cap (I + \overline{1/2}) = \emptyset$. Let $f(r, \theta) = (r',
\theta')$ be a homeomorphism of the cylinder that satisfies:

\begin{itemize}
\item[$(i)$] $\Delta_r = r' - r$ and $\Delta_{\theta} = \theta' - \theta$ only depend on $\theta$.
\item[$(ii)$] $\Delta_r = 4$ if $\theta \notin I$, $\Delta_r \in [-1, 0)$ if $\theta$ belongs to an
 interval $J \subset I$, $\bar{0} \in J$ and equals $-1$ if $\theta = \bar{0}$, see Figure \ref{fig:1}.
\item[$(iii)$] $\Delta_{\theta}$ is non-negative, $\Delta_{\theta} \le \dist(I, I+\overline{1/2})$
 and $\Delta_{\theta} = 0$ if and only if $\theta = \bar{0}$.
\end{itemize}

\begin{figure}[h!]
\begin{tabular}{c c c }
\includegraphics[scale = .8]{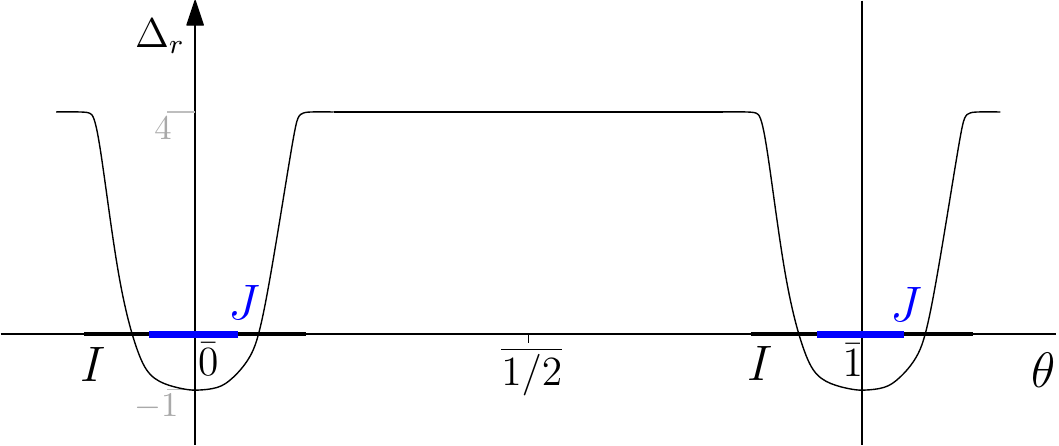} & \qquad & \includegraphics[scale = .8]{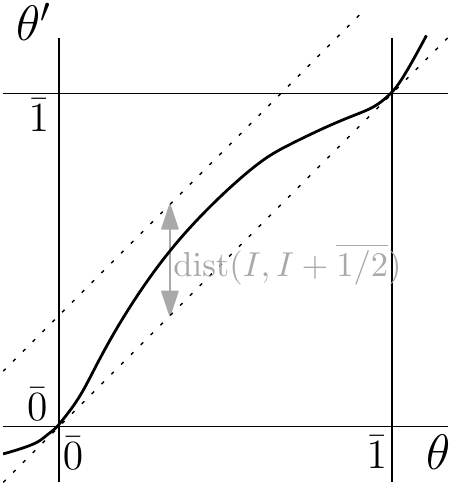}
\end{tabular}
\caption{Graphs of $\Delta_r$ (left) and $\theta'$ (right) as a
function of $\theta$ in the definition of $f.$}\label{fig:1}
\end{figure}

Property $(iii)$ controls the 1-dimensional dynamics in the angular
coordinate:  every orbit tends to $\theta = \bar{0}$. By $(ii)$ the
radial coordinate decreases indefinitely once the orbit remains
close to $\theta = \bar{0}$. From $(i)$ we deduce that $\Delta_r$ is
uniformly bounded and, as a consequence, $f$ extends to a planar
homeomorphism, which we will also denote by $f$, by imposing that
the origin is a fixed point $f(\textbf{0}) = \textbf{0}$.

The angular interval $J$ determines an infinite cone $\widehat J$ in which
the radial  coordinate of a point decreases after applying $f$.
Inside $\widehat J$ we find the trapping sector that has been
previously mentioned,
%it is determined by the positively invariant subinterval whose endpoint is $\bar{0},$
see Figure \ref{fig:2}.
Notice also that the speed of attraction ($\Delta_r \in [-1,0)$) in
the trapping sector is weaker than the speed of repulsion ($\Delta_r
= 4$) outside the cone $\widehat I$ determined by $I$.
%Indeed, most points in the plane are moved away from the origin after one iteration of $f$.

\begin{figure}[h!]
\includegraphics[scale = .8]{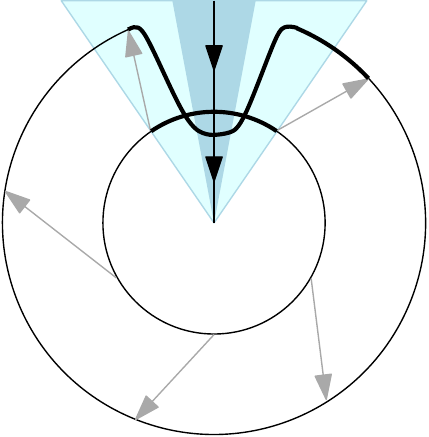}
\caption{Action of $f$ on a circle $r = c$ (inner circle) represented schematically by arrows, cones
$\widehat I$ (light) and $\widehat J$ (dark) are shadowed, the attracting $f$--invariant
ray ($\theta = \bar{0}$) is depicted vertical.}\label{fig:2}
\end{figure}

The map $g$ is merely a copy of $f$ shifted in the angular
coordinate.  Let $\tau(r, \theta) = (r, \theta + \overline{1/2})$ be
the half-turn rotation in the plane and let $g = \tau^{-1} \circ f
\circ \tau = \tau \circ f \circ \tau$. Being a conjugate of $f$, $g$
satisfies the same properties as $f$ if we replace $\theta$ by
$\theta + \overline{1/2}$ in the statements. The key observation is
that by the second item in $(iii)$,
$f(\widehat I) \cap \tau(\widehat I) = \{\textbf{0}\}$ and $g(\tau(\widehat I)) \cap \widehat I = \{\textbf{0}\}$.
This means that the radial coordinate cannote decrease under the action of $f$ and then immediately decrease under the action of $g$, or viceversa. In fact, by $(ii)$ the radial coordinate increases after applying $g \circ f$ or $f \circ g$.
%the image under $g$ of $\widehat J$ is disjoint to itself and, similarly, the image under $f$ of the attracting sector for $g$, $\tau(\widehat J)$ is disjoint to itself.
%This fact together with the condition on the speeds of convergence and divergence from the origin guarantee that the compositions $f \circ g$ and $g \circ f$ become repelling maps.

\smallskip

Let us finally prove that $f$ and $g$ exhibit the dynamical Parrondo's paradox. 

\smallskip

\textbf{The origin is a globally attracting fixed point for $f$ and $g$.}
%\medskip
Let $\{(r_n, \theta_n)\}_{n \ge 1}$ be the orbit under
$f$ of a  point $(r, \theta)$. We study separately the angular
coordinate $\{\theta_n\}$ because its evolution is independent from
the values of the radial coordinate. If $\theta = \bar{0}$ then
$\theta_n = \bar{0}$ for every $n \ge 1$, otherwise the sequence
$\{\theta_n\}$ is increasing and converges to the unique value
$\theta_0$ which is fixed under the 1-dimensional angular dynamics,
namely $\theta_0 = \bar{0}$. Thus, $\theta_n \to \bar{0}$ so
$\theta_n \in J$ for sufficiently large $n$, say $n \ge n_0$. This
implies that $r_{n+1} < r_n$ for every $n \ge n_0$ and,
additionally, that $r_{n+1} - r_n \to -1$ as $n \to +\infty$ because
the orbit converges towards the half-ray $\theta = \bar{0}$. As a
consequence, $r_n \to -\infty$ and the orbit tends to the
origin.

The same conclusion holds for $g$ because it is conjugate to $f$.

\smallskip
\textbf{The origin is a globally repelling fixed point for $f \circ g$ and $g \circ f$.}
%\medskip
Recall that $f \circ g = f \circ \tau \circ f \circ \tau$ is
conjugate to $g \circ f = \tau \circ f \circ \tau \circ f$ (notice
that $\tau^2 = \id$) so it is enough to prove the statement for the
latter composition. By $(ii)$ the radial coordinate of a point
outside $\widehat{I}$  increases by 4 under the action of $f$. Thus,
if $(r', \theta') = g \circ f (r, \theta)$ we have that $r' - r \ge
3$ if $\theta \notin I$. The same inequality is true in the case
$f(r, \theta)$ does not belong to $\tau(\widehat{I})$. Since
$\widehat{I} \cap f^{-1}(\tau(\widehat{I})) = \{\textbf{0}\}$ we conclude
that the radial coordinate of every point increases at least by 3
after applying $g \circ f$. Evidently, the origin is a global
repeller for $g \circ f$ and the proof of Theorem \ref{th:main1}
follows for $k=2.$

\subsection{Proof of Theorem \ref{th:main1} for $k>2$}

First, let us modify slightly the planar dynamics introduced in the previous subsection in order to make it symmetric with respect to the vertical axis. Define $h(r, \theta) = f(r, 2\theta)$ if $\theta \in [\bar{0},\overline{1/2}]$ and $h(r, \theta) = f(r, 1 - 2\theta)$ if $\theta \in [\overline{1/2}, \bar{1}]$. There are now two invariant rays for $h$, $\theta = \bar{0}$, which acts as a repeller in the dynamics in the angular coordinate, and $\theta = \overline{1/2}$, which acts as an attractor. The dynamics of $h$ within each half-plane reproduces the dynamics of $f$ except from the fact that both $h$-invariant rays correspond to the unique invariant ray $\{\theta = \bar{0}\}$ for $f$.

Now, it is straightforward to move into higher dimensions. Consider spherical coordinates $(r, \theta, \varphi_1, \ldots, \varphi_{k-2})$ in $\mathbb R^k$ and define a map $h_k \colon \mathbb R^k \to \mathbb R^k$ by the transformation that applies $h$ to the radial and polar coordinate, $(r', \theta') = h(r, \theta)$, and leaves the azimuthal coordinates unchanged, $(\varphi'_1, \ldots, \varphi'_{k-2}) = (\varphi_1, \ldots, \varphi_{k-2})$.
The dynamics of $h_k$ leaves invariant the two rays that form the vertical axis (north-south direction, suppose that north corresponds to $\theta = \bar{0}$). Points are attracted to the origin in those rays.
Moreover, the radial coordinate of any point increases significantly after applying $h_k$ unless the point belongs to a thin double cone $C$ around the axis (whose size can be traced back to the size of $I$). Nevertheless, since every orbit either belongs to the ray pointing to the north or eventually enters the cone around the ray that points to the south and remains in it, we conclude that the origin is a globally attracting fixed point for $h_k$.

An analogous construction for the second map as in the case $k=2$ works here as well. Let $\tau_k$ be a $90^o$-degree rotation in $\mathbb R^k$ and define $j_k = \tau_k^{-1} \circ h_k \circ \tau_k$. Note that

\centerline{
($\star$) \enskip \enskip
$h_k(C) \cap \tau_k(C) = \{\textbf{0}\}$ \enskip and \enskip $j_k(\tau_k(C)) \cap C = \{\textbf{0}\}$
}
\noindent
It is straightforward to check that the origin is a globally attracting fixed point for $j_k$ (again by conjugation) and that the radial coordinate of every point increases under the action of $j_k \circ h_k$ and $h_k \circ j_k$ (by ($\star$)) and the origin is a globally repelling fixed point for the composite maps.

\section{Iterated function system: proof of Theorem \ref{th:main2}}

The idea is to take as $f_0$ and $f_1$ a slight modification of the maps $f$ and $g$ defined in the proof of Theorem \ref{th:main1}. For the sake of clarity, we only discuss the case $k = 2$, the proof for $k > 2$ is a straightforward generalization of the argument using the maps $h_k$ and $j_k$.

Let us start with the proof of Theorem \ref{th:main2} for $k=2$.
% by introducing the maps $f_0 = f$ and $f_1 = g$.
%To prove Theorem \ref{th:main2} we start introducing both maps $f_0=f$ and $f_1=g.$
We need to slightly modify the definition of $f$ and $g$ in Subsection \ref{ss:fg} in order to increase the rate of radial repulsion far from the invariant rays to account for the effect of the probability $p$.
%They will be slight modifications of the $f$ and $g$ defined in previous section.
The only change in the definition of the new $f$, which we shall denote in the following by $f_0$,
is that we replace property $(ii)$ by
%will be that instead of property $(ii)$, the new $f$ satisfies property:
\begin{itemize}
\item[$(ii)'$] $\Delta_r = a-1,$ for some fixed $a>4,$ to be determined later, if $\theta \notin I$, $\Delta_r \in [-1, 0)$ if $\theta$ belongs to an
 interval $J \subset I$, $\bar{0} \in J$ and equals $-1$ if $\theta = \bar{0}$.
\end{itemize}
Then, the new $g$, which shall be henceforth denoted $f_1$, is constructed from the new $f$ as in the
previous section, $f_1 = \tau^{-1} \circ f_0 \circ \tau$. Notice that the original $f$ and $g$ considered in
Subsection \ref{ss:fg} correspond to $a=5.$ The value $a$ will be
fixed later, in terms of $p.$

Take an arbitrary point ${\bf x}$  in the plane, different from the
origin, and apply $f_0$ and $f_1$ randomly as in the statement, that is, apply
$f_0$ with probability $p$ and $f_1$ with probability $q = 1-p.$ We claim that the orbit of ${\bf x}$ is repelled from the origin almost surely,
that is, with probability 1.

The proof of the claim follows from two remarks. The first one
concerns the four  maps $f_0 \circ f_0, f_0 \circ f_1, f_1 \circ f_0, f_1 \circ
f_1$. Their radial coordinate change is bounded by:

\[\Delta_r^{f_0 \circ f_0} \ge -2,\quad \Delta_r^{f_0 \circ f_1} \ge
a-2,\quad  \Delta_r^{f_1 \circ f_0} \ge a-2\quad \mbox{and}\quad
\Delta_r^{f_1 \circ f_1} \ge -2.\] Moreover, the map $f_0\circ f_0$ appears
with probability $p^2,$ the map $f_1 \circ f_1$ with probability $q^2,$
while each of the maps $f_0\circ f_1$ and $f_1\circ f_0$ appears with
probability $pq.$ Let us start giving conditions on $a$ that imply
that the expected value of the change in radial coordinate is
positive. More precisely, if $\Delta^n$ denotes the random variable
that measures the minimum of the variation of radial coordinate
between a point (different from the origin) and its image under $F_{a_n,\cdots, a_1, a_0}$ we
have that
\begin{align*}
E[\Delta^{2m+2}] &\ge E[\Delta^{2m}] + p^2\min\Delta_r^{f_0 \circ f_0} +
pq\min\Delta_r^{f_0 \circ f_1}+  pq\min\Delta_r^{f_1 \circ f_0} +
q^2\min\Delta_r^{f_1 \circ f_1}\\&\ge E[\Delta^{2m}]+2(a-2)pq
-2(p^2+q^2)=E[\Delta^{2m}]+2\big(ap(1-p)-1\big).
\end{align*}
Hence, if we take any $a$ such that $ap(1-p)-1>0$
 we have that $E[\Delta^{2m+2}]\ge
E[\Delta^{2m}]+ K,$ for $K=2\big(ap(1-p)-1\big)>0,$ and as a
consequence we conclude that $E[\Delta^{2m}] \ge 2K m$, that is, the
expected value of $\Delta^n$ grows  linearly with $n$. This
computation shows that in average random iteration repels points
from the origin by increasing (linearly!) its radial coordinate. We
need to extend this conclusion to a subset of binary sequences of
full probability. Notice incidentally that $a>1/\big(p(1-p)\big)\ge 4.$

Given a binary sequence $(a_n)$ we can bound the value of $\Delta^{2m}$ in the following fashion
\[
\Delta^{2m} \ge \Delta_r^{f_{a_1}\circ f_{a_0}} +
\Delta_r^{f_{a_3}\circ f_{a_2}} + \ldots + \Delta_r^{f_{a_{2m}}\circ
f_{a_{2m-1}}} \ge (a-2)k_m - 2(m-k_m) = ak_m-2m,
\]
where $k_n$ denotes the number of maps among $f_{a_1} \circ f_{a_0},
f_{a_3} \circ f_{a_2},  \ldots, f_{a_{2n}} \circ f_{a_{2n-1}}$ that
are equal to $f_1 \circ f_0$ or $f_0 \circ f_1$. Notice that   $k_n$ is  the
sum of $n$ independent Bernoulli distributions
$\operatorname{B}(2pq),$ because $2pq$ is the probability of
appearance of  $f_1 \circ f_0$ or $f_0 \circ f_1$. Thus, if $\liminf_{n \to
+\infty} k_n/n = \ell > 2/a$, the asymptotic growth of $\Delta^{2n}$
is bounded from below by $(\ell - 2/a)n$. In particular,
$\Delta^{2n} \to +\infty$ so every point is repelled from the origin
by the iterated action of the maps $f_{a_n}, n \ge 1$.

It only remains to prove that the subset of $(a_n)$ such that
$\liminf k_n/n > 2/a$  has full probability. In fact, the  Strong
Law of Large Numbers (\cite{A,B}) gives much more: for a full
probability set, the previous $\liminf$ is indeed a limit and it  coincides
with the expected value of the random variable
$\operatorname{B}(2pq),$ that is $2pq.$ Hence, for a full
probability set of binary sequences,
\[
\ell=\liminf_{n \to +\infty}\frac {k_n}n=\lim_{n \to +\infty}\frac
{k_n}n=2pq=2p(1-p).
\]
For those sequences we have that,
\[
\ell-\frac 2a =2p(1-p)-\frac 2a=\frac{2\big(ap(1-p)-1\big)}{a}>0,
\]
as we wanted to prove, and the theorem follows.

\subsection*{Acknowledgements} This work has received funding from the Ministerio de Ciencia e Innovaci\'{o}n
(MTM2016-77278-P FEDER, PGC2018-098321-B-I00 and PID2019-104658GB-I00 grants), the Ag\`{e}ncia
de Gesti\'{o} d'Ajuts Universitaris i de Recerca (2017 SGR 1617 grant).


\begin{thebibliography}{12}

\bibitem{A}  R. B. Ash. Real analysis and probability. Probability and Mathematical Statistics, No. 11. Academic Press, New York-London, 1972. %xv+476 pp.


\bibitem{B} P. Billingsley. Probability and measure. Third edition.
Wiley Series in Probability and Mathematical Statistics. A
Wiley-Interscience Publication. John Wiley \& Sons, Inc., New York,
1995.% xiv+593 pp.



    \bibitem{BTT} V. D.~Blondel, J.~Theys, J. N.~Tsitsiklis.\textsl{When is a pair of matrices
        stable?}. In: V. D.~Blondel, A. Megretski (eds.). Unsolved problems
    in Mathematical Systems and Control Theory. Princeton Univ. Press,
    NJ 2004.



    \bibitem{CLP} J. S.~C\'anovas, A.~Linero, D.~Peralta-Salas. \textsl{Dynamic Parrondo's paradox.}
    Physica D 218 (2006) 177--184.



    \bibitem{CGM18} A.~Cima, A.~Gasull, V.~Ma\~{n}osa.
    Parrondo's dynamic paradox for the stability of non-hyperbolic fixed
    points.  Discrete Contin. Dyn. Syst.  38 (2018), 889--904.



    \bibitem{ES1} S.~Elaydi, R. J.~Sacker.
    \textsl{Global stability of periodic orbits of non-autonomous
        difference equations and population biology.} J. Differential
    Equations 208 (2005), 258--273.


    \bibitem{ES2} S.~Elaydi, R. J.~Sacker. \textsl{Periodic
        difference equations, population biology and the Cushing-Henson
        conjectures.} Math. Biosci. 201 (2006), 195--207.

    \bibitem{FS} J. E.~Franke, J. F.~Selgrade. \textsl{Attractors for discrete periodic
        dynamical systems.} J. Math. Anal. Appl. 286 (2003), 64--79.

    \bibitem{HA}  G. P. Harmer and D. Abbott. \textsl{Losing strategies can win by Parrondo's paradox.}
    Nature (London), Vol. 402, No. 6764 (1999) p. 864.


    \bibitem{J} R.~Jungers. The Joint Spectral Radius. Spinger,
    Berlin 2009.



    \bibitem{Par} J. M. R. Parrondo. \textsl{How to cheat a bad
        mathematician.} in EEC HC\&M Network on Complexity and Chaos
    (\#ERBCHRX-CT940546), ISI, Torino, Italy (1996), Unpublished.

%{\color{blue} \bibitem{Sed} H. Sedaghat. \textsl{The Impossibility of Unstable, Globally Attracting Fixed Points for Continuous Mappings of the Line.} Amer. Math. Monthly 104 (1997), 356--358.}


    \bibitem{SR1} J. F.~Selgrade, J. H.~Roberds. \textsl{On the structure of
        attractors for discrete, periodically forced systems with
        applications to population models.} Physica D 158 (2001), 69--82.

    \bibitem{SR2} J. F.~Selgrade, J. H.~Roberds. \textsl{Global attractors for a
        discrete selection model with periodic immigration.} J. Difference Equations and Appl. 13 (2007), 275--287.

%{\color{blue} \bibitem{SMR} A. N. Sharkovsky, Yu. L. Maistrenko, and E. Yu. Romanenko. Difference Equations and Their Applications. Kluwer Academic Publ., Dordrecht, 1993. }

\end{thebibliography}
\end{document}